# SUBSPACE ESTIMATION AND PREDICTION METHODS FOR HIDDEN MARKOV MODELS

By Sofia Andersson and Tobias Rydén[1]

*AstraZeneca R&D and Lund University*

Hidden Markov models (HMMs) are probabilistic functions of finite Markov chains, or, put in other words, state space models with finite state space. In this paper, we examine subspace estimation methods for HMMs whose output lies a finite set as well. In particular, we study the geometric structure arising from the nonminimality of the linear state space representation of HMMs, and consistency of a subspace algorithm arising from a certain factorization of the singular value decomposition of the estimated linear prediction matrix. For this algorithm, we show that the estimates of the transition and emission probability matrices are consistent up to a similarity transformation, and that the $m$-step linear predictor computed from the estimated system matrices is consistent, i.e., converges to the true optimal linear $m$-step predictor.

**1. Introduction.** A hidden Markov model (HMM) is a time series model for which the observable output is governed by a nonobservable finite-state Markov chain. HMMs have been applied in a wide range of areas, including telecommunications, speech recognition, bioinformatics and econometrics. See the monographs [8, 12, 19, 22] as well as the somewhat old but still informative tutorial [28] for further reading and references.

Estimation of HMMs is, in most cases, carried out using maximum likelihood (ML), and it is known that the MLE enjoys the standard favorable statistical properties: consistency, asymptotical normality and asymptotic efficiency [5, 7, 21]. However, from a computational point of view, the MLE is less tractable; the log-likelihood function is highly nonlinear and typically

Received October 2008; revised April 2009.

[1]Supported by a grant from the Swedish Natural Science Research Council.

*AMS 2000 subject classifications.* Primary 62M09; secondary 62M10, 62M20, 93B15, 93B30.

*Key words and phrases.* Hidden Markov model, linear innovation representation, prediction error representation, subspace estimation, consistency.

severely multi-modal. Thus, iterative algorithms like the EM algorithm, often referred to as the Baum–Welch algorithm in the context of HMMs [6], are needed to maximize the log-likelihood, and there is no guarantee that such algorithms will find the actual global maximum—the MLE.

Subspace methods have proved very successful for estimating linear Gaussian state-space models. These methods are noniterative and based on matrix computations that are numerically stable. An HMM can be formulated as a linear state-space model, and in [16], subspace estimation was applied to a simple HMM; no further theoretical basis for the validity of the procedure was provided, but encouraging numerical results were reported. The purpose of the present paper is to give a theoretical foundation for the application of subspace methods to HMMs. Our main results are two consistency theorems. The first result states that the estimated transition and emission probability matrices are consistent, up to a similarity transformation (change of basis). The second result states that the optimal linear $m$-step predictor computed from the estimated matrices, is consistent in the sense that it converges to the corresponding optimal predictor given the true, unknown, system matrices. This second result can also be viewed as HMMs, estimated with subspace methods, being a general tool for very quickly (complexity discussed in next paragraph) building predictors for discrete time series (or, more precisely, time series taking values in a finite set). In this respect, the subspace approach is analogous to building predictors with finite-dimensional state-space models (or, equivalently, ARMA models) for real-valued time series.

Difficulties remain however, in particular the so-called positive realization problem, which here amounts to determining a representation of an HMM second-moment function in terms of transition probability matrices which are positive and with elements between zero and one. Significant progress has been made on this subject [1, 30], and there exist locally convergent algorithms to find positive realizations for related problems (cf. Section 9), but it is fair to say that the problem is not settled. Another drawback of subspace methods is that they do not allow for including structural constraints on the transition and/or emission probability matrix. For instance, it may be that some elements of these matrices are identically zero (by modeling assumptions), or that several probabilities are jointly given by a common lower-dimensional parameter. Such structures are usually easily incorporated by, for example, the EM algorithm, but cannot be handled with subspace methods. See Section 9 for a remark on zero elements. The advantage of subspace methods is rather on the computational side. With ML, one evaluation of the log-likelihood (amounting to running the forward pass of the forward–backward algorithm), or one iteration of EM (which is dominated by the complexity of the E-step) has complexity $O(Tn^2)$ with $T$ being the number of observations and $n$ being the cardinality of the state space ([8], Section 5.1.1) and many evaluations or iterations need to be carried out.

A subspace estimation algorithm first builds covariance matrices of size $(k\ell)^2$ for some $k$ from $T$ observations, where $\ell$ is the cardinality of the set of observations, and then computes the singular value decomposition of a matrix of the same size; the complexity is thus $O(T(k\ell)^2) + O((k\ell)^3)$ (cf. [13], Section 5.4.5). Typically, $k$ is chosen of the order $\log T$, and the complexity is then slightly higher than linear in $T$. In practice, a subspace method is much faster than a likelihood-based method however, because it is noniterative and all steps are carried out just once. This difference becomes particularly pronounced when $n$ and/or $\ell$ is large, since likelihood-based methods then in general suffer the most from multi-modal log-likelihood surfaces, require many iterations to converge, and are prone to converge to local maxima.

We do not give a general introduction to subspace methods in this paper, but refer the reader to, for example, [18, 26]. A problem we do not consider is order estimation, that is, estimating the number of states of the hidden Markov chain. It has been suggested in [27], Chapter 8, to use the singular value decomposition of the subspace algorithm for this purpose, by estimating the number of nonzero singular values. It is, however, not obvious that this approach immediately carries over from linear Gaussian models to HMMs. In the present paper, we leave this issue alone assuming that the model order is known beforehand or estimated otherwise. Another topic that we will not elaborate on is the choice of truncation index, or time horizon, $k$, of the subspace algorithm. Although our consistency result gives guidelines on how large, asymptotically, $k$ is allowed to be, it gives no information on what a good choice for $k$ for a finite set of data would be. It has been proposed to select $k$ as to minimize certain criteria—see [4], Chapter 5, and [27], Chapter 8—but again we do not consider this problem further.

The paper is organized as follows. Section 2 sets the notation for the HMM and some matrices used by the algorithm, and Section 3 describes various state-space representations of the model. Section 4 discusses linear prediction of future observations from past ones, while Section 5 summarizes the algorithm. Section 6 explores structural properties of the system related to its nonminimality, and Section 7 contains consistency results. The paper is concluded with a small simulation study in Section 8, and a discussion of the theoretical results in Section 9.

**2. Notation.** We consider an HMM with finite state space and observations in a finite alphabet. More specifically, the model comprises a nonobservable (hidden) Markov chain $\{x_t^M\}_{t=-\infty}^{\infty}$ (the state process) and an observable discrete-time process $\{y_t\}_{t=-\infty}^{\infty}$ such that (i) given $\{x_t^M\}$, $\{y_t\}$ is a sequence of conditionally independent random variables, and (ii) the conditional distribution of $y_{t+1}$ depends on $x_t^M$ only. The state space will be denoted by $\{e_1, \ldots, e_n\}$ (its size is $n$) and the output alphabet will be denoted by $\{d_1, \ldots, d_\ell\}$ (its size is $\ell$). Assuming that $\{x_t^M\}$ is stationary, the model is completely characterized by the $n \times n$ transition probability matrix $A$ and the $\ell \times n$ emission probability matrix $C$, with elements

$$a_{ij} = \mathbb{P}(x_{t+1}^M = e_i | x_t^M = e_j), \qquad c_{ij} = \mathbb{P}(y_{t+1} = d_i | x_t^M = e_j);$$

here $\mathbb{P}$ denotes probability. This notation is somewhat nonstandard in that usually $a_{ij}$ is the probability of moving *from $e_i$ to $e_j$* (but here the opposite) and $x_t^M$ governs the conditional distribution of $y_t$ (but here $y_{t+1}$). As will become clear below, however, the present notation leads to convenient expressions in what follows.

Without loss of generality, we will let $e_i$ be the $i$th coordinate (column) vector in $\mathbf{R}^n$. $e_i \in \mathbf{R}^n$ thus consists of zeros, except for its $i$th element which is unity. Similarly, $d_i$ will be taken as the $i$th unit vector in $\mathbf{R}^\ell$. In fact, this choice is crucial in order to represent the HMM as a linear state-space model. Both $A$ and $C$ have entries in $[0, 1]$ and column sums equal to unity; we say that they are column stochastic. These constraints can be written $\mathbf{1}_n^\top A = \mathbf{1}_n^\top$ and $\mathbf{1}_\ell^\top C = \mathbf{1}_n^\top$, where $\mathbf{1}_n$ is a length $n$ column vector of ones and superindex $\top$ denotes matrix transposition. Also note that $\mathbf{1}_n^\top x_t^M = 1$ and $\mathbf{1}_\ell^\top y_t = 1$. We will write $\mathbf{S}_n$ for the linear subspace of $\mathbf{R}^n$ spanned by $\mathbf{1}_n$, and $\mathbf{S}_n^\perp$ for its orthogonal complement. We will write $\mathbb{E}$ for expectation, and $\mathbf{1}_n$ (as above) and $\mathbf{0}_n$ for a length $n$ column vector of ones or zeros, respectively. The following assumption will be imposed.

CONDITION A. *$A$ is irreducible and aperiodic (i.e., ergodic).*

Under this assumption, $A$ admits a unique stationary distribution $\pi$ satisfying $A\pi = \pi$. Note that $\pi = \mathbb{E}[x_t^M]$ and, similarly, $\mathbb{E}[y_t] = C\pi$. Moreover, $A$ has one eigenvalue on the complex unit circle (namely 1), while its remaining eigenvalues lie inside it. This implies that $A^t \to \pi \mathbf{1}_n^\top$ as $t \to \infty$.

**3. State-space representations.** In order to apply subspace methods, we need to formulate the HMM in state-space form. In [12], Chapter 2, it was shown that the model may be written as

$$(3.1) \qquad x_{t+1}^M = A x_t^M + \xi_{t+1},$$

$$(3.2) \qquad y_{t+1} = C x_t^M + \eta_{t+1},$$

where $\{\xi_t\}$ and $\{\eta_t\}$ are independent sequences of martingale increments. Hence, both processes contain uncorrelated elements, but neither process has elements that are independent or Gaussian. Iterating (3.1) leads to $x_t^M = \sum_0^\infty A^j \xi_{t-j} + \pi$, where the sum is well defined because $\mathbf{1}_n^\top \xi_t = 0$. This system also admits a prediction error (or, innovation) representation [3],

$$(3.3) \qquad x_{t+1} = A x_t + K \varepsilon_{t+1},$$

$$(3.4) \qquad y_{t+1} = C x_t + \varepsilon_{t+1}.$$

Here, $x_t$ is the optimal linear predictor of $x_t^M$ given all past observable information $\{y_s\}_{s=-\infty}^t$. Note that $\mathbf{1}_n^\top x_t = 1$. Moreover, $\varepsilon_t$ is the one-step-ahead linear prediction error of $y_t$ given $\{y_s\}_{s=-\infty}^{t-1}$. Hence, $\mathbb{E}\varepsilon_t = 0$ and $\{\varepsilon_t\}$ is an uncorrelated sequence. The $n \times \ell$ matrix $K$, the Kalman gain, solves a "singular" Riccati equation and satisfies $\mathbf{1}_n^\top K = 0$. Finally, define the centered processes $\overline{x}_t = x_t - \pi$ and $\overline{y}_t^M = y_t^M - C\pi$, both of which have zero mean, as well as the matrix $\overline{A} = A - \pi \mathbf{1}_n^\top$. Subtraction by $\pi \mathbf{1}_n^\top$ moves the eigenvalue 1 of $A$ to 0, but otherwise leaves the eigenvalues and -vectors of $A$ alone. Thus, $\overline{A}$ is a stable matrix, that is, its spectral radius is less than one. Using $A\pi = \pi$ and $\mathbf{1}_n^\top \overline{x}_t = 0$, we find

$$(3.5) \qquad \overline{x}_{t+1} = A\overline{x}_t + K\varepsilon_{t+1} = \overline{A}\overline{x}_t + K\varepsilon_{t+1},$$

$$(3.6) \qquad \overline{y}_{t+1} = C\overline{x}_t + \varepsilon_{t+1}.$$

Since $\mathbf{1}_n^\top x_t = 1$, the $x$-process lives in an $(n-1)$-dimensional affine subspace of $\mathbf{R}^n$. A major idea of the present paper is to restructure this process into one component $\overline{x}$ living in the $(n-1)$-dimensional linear space $\mathbf{S}_n^\perp$, and one component living in the null-dimensional affine space $\{1\}$. Thereby, we need to estimate the $(n-1)$-dimensional component in the linear space only, whereas the remaining affine component is for free.

We now introduce the notation

$$y_t^+ = (y_{t+1}^\top, y_{t+2}^\top, \ldots)^\top, \qquad y_t^+(k) = (y_{t+1}^\top, y_{t+2}^\top, \ldots, y_{t+k}^\top)^\top,$$

$$y_t^- = (y_t^\top, y_{t-1}^\top, \ldots)^\top, \qquad y_t^-(k) = (y_t^\top, y_{t-1}^\top, \ldots, y_{t-k+1}^\top)^\top.$$

These quantities are all column vectors. Vectors $\varepsilon_t^+$, $\varepsilon_t^-$, etc., are defined in a corresponding manner, as are the vectors $\overline{y}_t^+$, $\overline{y}_t^-$, etc., of centred observations (subtracting $C\pi$). Further, we define the covariance matrices

$$\mathcal{H} = \mathbb{E}[\overline{y}_t^+ \overline{y}_t^{-\top}], \qquad \mathcal{H}_k = \mathbb{E}[\overline{y}_t^+(k) \overline{y}_t^-(k)^\top],$$

$$\Gamma = \mathbb{E}[\overline{y}_t^- \overline{y}_t^{-\top}], \qquad \Gamma_k = \mathbb{E}[\overline{y}_t^-(k) \overline{y}_t^-(k)^\top].$$

Here, $\mathcal{H}_k$ contains the first $k$ (size $\ell$) block rows and block columns of $\mathcal{H}$, and similarly for $\Gamma_k$. The entries of these matrices can be expressed using

$$\mathbb{E}[\overline{y}_t \overline{y}_{t+j}^\top] = CA^{(-j\vee 0)} SA^{(j\vee 0)\top} C^\top + \delta_{j0} R = C\overline{A}^{(-j\vee 0)} S\overline{A}^{(j\vee 0)\top} C^\top + \delta_{j0} R,$$

where $\delta_{j0}$ is Kronecker's delta and

$$S = \mathbb{E}[(x_t^M - \pi)(x_t^M - \pi)^\top] = \operatorname{diag}(\pi) - \pi\pi^\top,$$

$$R = \mathbb{E}[\eta_t \eta_t^\top] = \operatorname{diag}(C\pi) - C\operatorname{diag}(\pi)C^\top.$$

The equalities above follow from (3.1) and (3.2); alternatively, the second equality can be derived by observing (e.g., using induction) that $A^j = \overline{A}^j + \pi \mathbf{1}_n^\top$ for any $j \geq 0$, and that $\mathbf{1}_n^\top S = (S\mathbf{1}_n)^\top = 0$.

**4. Projecting future observations on past ones.** Equations (3.5) and (3.6) yield

$$\overline{y}_{t+m} = C\overline{A}^{m-1}\overline{x}_t + \sum_{j=1}^{m-1} C\overline{A}^{j-1}K\varepsilon_{t+m-j} + \varepsilon_{t+m}. \tag{4.1}$$

Defining $\mathcal{O} = (C^\top \ \overline{A}^\top C^\top \ (\overline{A}^2)^\top C^\top \ \ldots)^\top$, we can write this in matrix form as $\overline{y}_t^+ = \mathcal{O}\overline{x}_t + \mathcal{M}\varepsilon_t^+$ for some matrix $\mathcal{M}$. Now, write $\mathbb{H}_{\overline{y}_t^-}$ for the Hilbert space spanned by the elements of $\overline{y}_t^-$ and $\mathcal{P}_{\overline{y}_t^-}$ for the projection operator onto this space. Since $\overline{x}_t \in \mathbb{H}_{\overline{y}_t^-}$ by definition and $\varepsilon_t^+$ is orthogonal to $\overline{y}_t^-$, we obtain $\mathcal{P}_{\overline{y}_t^-}\overline{y}_t^+ = \mathcal{O}\overline{x}_t$. Furthermore, it follows from [3] [see (3.5) and (3.6)] that

$$\overline{x}_t = \sum_{j=0}^{\infty}(A-KC)^j K\overline{y}_{t-j} = \mathcal{K}\overline{y}_t^-, \tag{4.2}$$

where $\mathcal{K} = (K, (A-KC)K, (A-KC)^2 K, \ldots)$. Hence, $\mathcal{P}_{\overline{y}_t^-}\overline{y}_t^+ = \mathcal{O}\mathcal{K}\overline{y}_t^- = \beta\overline{y}_t^-$, with, obviously, $\beta = \mathcal{O}\mathcal{K}$. Thus, $\beta$ is the matrix that describes the optimal linear prediction of $\overline{y}_t^+$ from $\overline{y}_t^-$. We note that $\beta = \mathcal{H}\Gamma^\dagger$, where superscript $\dagger$ denotes the Moore–Penrose pseudo-inverse; that $\Gamma^\dagger$ is well defined is proved in Section 6.

Now, assume that we want to predict in a finite horizon context, that is, we want to predict linearly $\overline{y}_t^+(k)$ from $\overline{y}_t^-(k)$ for some $k$. The optimal weight matrix $\beta_k$ is then $\beta_k = \mathcal{H}_k\Gamma_k^\dagger$. It is indeed necessary to use this generalized inverse since $\Gamma_k$ and other matrices are singular as a result of the nonminimality of the system, to be discussed in detail in Section 6. Note that $\beta_k$ is not an approximation of $\beta$, as the latter matrix is infinite whereas $\beta_k$ is $k\ell \times k\ell$, but we will show below that the predictions obtained using these two matrices become identical as $k \to \infty$.

The following assumption ensures that the matrix $\beta$ has maximal rank, apart from a rank loss of one due to the nonminimality of the system.

CONDITION B. The matrix $\mathcal{O}$ has rank $n$, and the matrix $\mathcal{K}$ has rank $n-1$.

This condition thus states that $\mathcal{O}$ is full rank, and that the relation $\mathbf{1}_n^\top K = 0$ is the only rank deficiency of $\mathcal{K}$. It ensures that $\beta = \mathcal{O}\mathcal{K}$ has rank $n-1$. In terms of linear systems theory, the first part of the assumption is equivalent to saying that the *observability matrix* formed by the pair $(A,C)$, that is, the $n\ell \times n$ matrix formed by the first $n$ block rows $CA^{i-1}$, has rank $n-1$ (cf. [29], Theorem 9.11); this is a straightforward application of the Cayley–Hamilton theorem. The second part amounts to saying that the *controllability matrix* formed by the pair $(A - KC, K)$, that is, the $n \times n\ell$ matrix formed by the first $n$ block columns $(A - KC)^{i-1}K$ of $\mathcal{K}$, has rank $n - 1$ (cf. [29], Theorem 9.5). This is equivalent to assuming that the controllability matrix formed by the pair $(A, K)$ has rank $n - 1$ (see [29], Theorem 18.16).

**5. The algorithm.** Assume that we have observed $\{y_t\}_{t=1}^T$ and define, for some $k \geq 1$, the data Hankel matrices

$$Y^+ = \begin{pmatrix} y_1 & y_2 & \cdots & y_T \\ y_2 & y_3 & \cdots & y_{T+1} \\ \vdots & \vdots & & \vdots \\ y_k & y_{k+1} & \cdots & y_{T+k} \end{pmatrix}, \qquad Y^- = \begin{pmatrix} y_0 & y_1 & \cdots & y_{T-1} \\ y_{-1} & y_0 & \cdots & y_{T-2} \\ \vdots & \vdots & & \vdots \\ y_{-k+1} & y_{-k+2} & \cdots & y_{T-k} \end{pmatrix},$$

where entries $y_s$ with $s \leq 0$ or $s > T$ are set to $d_i$ for arbitrary $i$. Empirically centred variables are defined by letting $m_{Y^+}$ and $m_{Y^-}$ be the row-wise means of $Y^+$ and $Y^-$, respectively, and putting $\overline{Y}^+ = Y^+ - m_{Y^+}\mathbf{1}_T^\top$ and $\overline{Y}^- = Y^- - m_{Y^-}\mathbf{1}_T^\top$. Then $\widehat{\mathcal{H}}_k = T^{-1}\overline{Y}^+\overline{Y}^{-\top}$ and $\widehat{\Gamma}_k = T^{-1}\overline{Y}^-\overline{Y}^{-\top}$ are estimates of $\Gamma_k$ and $\mathcal{H}_k$, respectively, and the standard least-squares estimate of $\beta_k$ is

$$\widehat{\beta}_k = (\overline{Y}^+\overline{Y}^{-\top})(\overline{Y}^-\overline{Y}^{-\top})^\dagger = \widehat{\mathcal{H}}_k\widehat{\Gamma}_k^\dagger.$$

Before proceeding, we need to reduce the rank of $\widehat{\beta}_k$, which is $k(\ell - 1)$, to that of the true $\beta$, which is $n - 1$. This is typically done by means of a singular value decomposition (SVD). Let

$$\widehat{\beta}_k = \widehat{U}\widehat{\Lambda}\widehat{V}^\top = (\widehat{U}_1, \widehat{U}_2)\begin{pmatrix} \widehat{\Lambda}_{11} & 0 \\ 0 & \widehat{\Lambda}_{22} \end{pmatrix}(\widehat{V}_1, \widehat{V}_2)^\top = \widehat{U}_1\widehat{\Lambda}_{11}\widehat{V}_1^\top + \widehat{U}_2\widehat{\Lambda}_{22}\widehat{V}_2^\top,$$

where $\widehat{\Lambda}_{11}$ and $\widehat{\Lambda}_{22}$ are $(n-1) \times (n-1)$ and $(k\ell - n + 1) \times (k\ell - n + 1)$ matrices, respectively be the SVD of $\widehat{\beta}_k$. We define estimates $\widehat{\mathcal{O}}_k$ and $\widehat{\mathcal{K}}_k$, of size $k\ell \times (n-1)$ and $(n-1) \times k\ell$, respectively to satisfy $\widehat{\mathcal{O}}_k\widehat{\mathcal{K}}_k = \widehat{U}_1\widehat{\Lambda}_{11}\widehat{V}_1^\top$, and our particular factorization of this equality will be

$$\widehat{\mathcal{O}}_k = \widehat{U}_1, \qquad \widehat{\mathcal{K}}_k = \widehat{\Lambda}_{11}\widehat{V}_1^\top.$$

Several other factorizations are more common in the literature, but this one is relatively simple to analyze as will be illustrated below. Our framework is however applicable to the analysis of other methods as well, such as for instance the canonical variate analysis (CVA) method (see [10]). Note that $\widehat{\beta}_k$ is an estimate before model order reduction and that $\widehat{\mathcal{O}}_k\widehat{\mathcal{K}}_k$ is an estimate after reduction. Thus, while $\beta = \mathcal{O}\mathcal{K}$, $\widehat{\beta}_k$ and $\widehat{\mathcal{O}}_k\widehat{\mathcal{K}}_k$ are not equal.

Another aspect of this factorization is that while the product $\widehat{\mathcal{O}}_k\widehat{\mathcal{K}}_k$ approximates a finite block of the infinite matrix $\mathcal{O}\mathcal{K}$, the factors $\widehat{\mathcal{O}}_k$ and $\widehat{\mathcal{K}}_k$ do not approximate $\mathcal{O}$ and $\mathcal{K}$, respectively in the same sense. As a result, the estimated states $\widehat{\overline{x}}_t$ below differ from the $\overline{x}_t$ by a change of basis, and the estimated system matrices $\widehat{A}$, etc., will differ from their true values by the same change of basis. Thus, the matrix $\widehat{A}$ will in general not have values in $[0,1]$, or column sums being unity; to find such a representation of the system requires computing a corresponding similarity transformation. An alternative way to view this problem is by noting that in (3.3) and (3.4), we can multiply the states $x_t$ and innovations $\varepsilon_t$ by any nonsingular matrix $S$, without changing the system if the matrices $A$, $C$ and $K$ are transformed accordingly. Predictive distributions computed using the estimated system matrices are consistent in the usual sense, however, without a change of basis.

Returning to the algorithm, an estimate $\widehat{\overline{X}} = (\widehat{\overline{x}}_0, \ldots, \widehat{\overline{x}}_{T-1})$ of the unobserved centred (predicted) state sequence $(\overline{x}_0, \overline{x}_2, \ldots, \overline{x}_{T-1})$ is constructed as $\widehat{\overline{X}} = \widehat{\mathcal{K}}_k(Y^- - (\mathbf{1}_k \otimes m_Y)\mathbf{1}_T^\top)$ [cf. (4.2)] where $\otimes$ denotes Kronecker product and $m_Y = T^{-1}\sum_1^T y_t$. Subtracting a common mean (of all observed $y_t$) makes the estimator slightly easier to analyze. Approximating $\overline{x}_t$ by $\widehat{\overline{x}}_t$ has two sides: first $\mathcal{K}$ is replaced by a matrix $\mathcal{K}_k$ comprising the first $k$ blocks of $\mathcal{K}$, and then $\mathcal{K}_k$ is replaced by an estimate. The estimate $\widehat{\overline{x}}$ is expressed in a different basis than is $\overline{x}$; note, however, that $\widehat{\overline{x}}$ and $\overline{x}$ lie in linear spaces $\mathbf{R}^{n-1}$ and $\mathbf{S}_n^\perp$, respectively, of equal dimensionality. An estimate of $x_t$ is obtained by adjoining the affine component 1, that is, we set $\widehat{x}_t = [\widehat{\overline{x}}_t^\top, 1]^\top$ and define $\widehat{X}$ as the corresponding $n \times T$ matrix. Again, $\widehat{x}_t$ differs from $x_t$ by a change of basis.

An estimate $\widehat{X}_1$ of $(x_1, x_2, \ldots, x_T)$ is obtained by dropping the first column of $\widehat{X}$ and then replicating its last column once. Our approximation of the system (3.3)–(3.4), up to a change of basis, is then

$$\widehat{X}_1 = A\widehat{X} + K\mathcal{E}, \qquad Y^+ = C\widehat{X} + \mathcal{E}.$$

Estimating $A$, $C$ and $K$ from these equations using linear regression gives

$$\widehat{A} = \widehat{X}_1 \widehat{X}^\top (\widehat{X}\widehat{X}^\top)^\dagger, \qquad \widehat{C} = Y^+ \widehat{X}^\top (\widehat{X}\widehat{X}^\top)^\dagger, \qquad \widehat{K} = \widehat{X}_1 \widehat{\mathcal{E}}^\top (\widehat{\mathcal{E}}\widehat{\mathcal{E}}^\top)^\dagger,$$
(5.1)
where $\widehat{\mathcal{E}}$ is the residual vector $\widehat{\mathcal{E}} = Y^+ - \widehat{C}\widehat{X}$; note that the regressions to find $\widehat{A}$ and $\widehat{K}$, respectively can be done separately since $\widehat{\mathcal{E}}$ is, by its definition as residuals, uncorrelated with $\widehat{X}$ (cf. [10], page 1867).

We remark that as the $n$th component of $\widehat{x}_t$ equals 1, it is clear that the estimated regression coefficients for this component, that is, the $n$th row of $\widehat{A}$, will be $[\mathbf{0}_{n-1}^\top, 1]$. Thus, $[\mathbf{0}_{n-1}^\top, 1]$ is a left eigenvector of $\widehat{A}$ with eigenvalue 1, so that $\widehat{A}$, like $A$, has one eigenvalue that is 1.

**6. Structural properties.** This section explores some structural properties of the state-space form of the HMM. The crucial difference between a "standard" state-space model and the model (3.1)–(3.2) is that the latter is not minimal. As noted above, $x_t^M$ lies in an $(n-1)$-dimensional affine subspace of $\mathbf{R}^n$ and, likewise, $y_t$ lies in an $(\ell-1)$-dimensional affine subspace. In addition, $A$ has a nonstable eigenvalue 1. Recall that $K$ satisfies $\mathbf{1}_n^\top K = 0$; also, $A - KC$ has a single eigenvalue 1 on the unit circle and its other eigenvalues inside the unit circle [3]. Obviously, $\mathbf{1}_n^\top(A - KC) = \mathbf{1}_n^T$, so that $\mathbf{1}_n^\top$ is the left eigenvector corresponding to the eigenvalue 1. Letting $\gamma$ be the corresponding right eigenvector, normalized so that $\mathbf{1}_n^\top \gamma = 1$, that is, $(A - KC)\gamma = \gamma$, we can define the matrix $J = (A - KC) - \gamma \mathbf{1}_n^\top$. This matrix plays a role similar to that of $\overline{A}$. In particular, $J$ is stable and $(A - KC)^j = J^j + \gamma \mathbf{1}_n^\top$ for any $j > 0$. Thus, the blocks $(A - KC)^j K$ that build $\mathcal{K}$ equal $J^j K$, and hence decay geometrically fast as $j \to \infty$.

Now define an $n \times (n-1)$ matrix $U_n = [u_1, u_2, \ldots, u_{n-1}]$ with orthonormal columns that span $\mathbf{S}_n^\perp$. The vector $n^{-1/2}\mathbf{1}_n$ has unit length and so $\Pi_{\mathbf{S}_n} = n^{-1}\mathbf{1}_n\mathbf{1}_n^\top$ is the matrix projecting a vector in $\mathbf{R}^n$ on $\mathbf{S}_n$ and $\Pi_{\mathbf{S}_n^\perp} = I_n - n^{-1}\mathbf{1}_n\mathbf{1}_n^\top$ is the matrix projecting on $\mathbf{S}_n^\perp$. Alternatively, $\Pi_{\mathbf{S}_n^\perp} = U_n U_n^\top$. The subspace $\mathbf{S}_n^\perp$ is isomorphic to $\mathbf{R}^{n-1}$, and this isomorphism is represented by the mappings $U_n : \mathbf{R}^{n-1} \to \mathbf{S}_n^\perp$ and $U_n^\top : \mathbf{S}_n^\perp \to \mathbf{R}^{n-1}$.

If $x \in \mathbf{S}_n^\perp$, then $\mathbf{1}_n^\top Ax = \mathbf{1}_n^\top x = 0$, that is, $Ax \in \mathbf{S}_n^\perp$. Thus, $A$ is a linear map on $\mathbf{S}_n^\perp$. Likewise, $C$ is a linear map from $\mathbf{S}_n^\perp$ to $\mathbf{S}_\ell^\perp$ and the covariance matrices $S$ and $R$ are linear operators on $\mathbf{S}_n^\perp$ and $\mathbf{S}_\ell^\perp$, respectively, since $\mathbf{1}_n^\top S = 0$ and $\mathbf{1}_\ell^\top R = 0$. The matrices $\widetilde{A} = U_n^\top A U_n = U_n^\top \overline{A} U_n$, $\widetilde{S} = U_n^\top S U_n$, $\widetilde{R} = U_\ell^\top R U_\ell$ and $\widetilde{C} = U_\ell^\top C U_n$ are the corresponding linear mappings operating between the isomorphic spaces $\mathbf{R}^{n-1}$ and $\mathbf{R}^{\ell-1}$. In particular, $\widetilde{A}$ is an $(n-1) \times (n-1)$ stable matrix. Define for any $k \geq 1$ the $k\ell \times k(\ell-1)$ block diagonal matrix

$$U_\ell^{[k]} = \begin{pmatrix} U_\ell & 0_{\ell \times \ell-1} & \cdots & 0_{\ell \times \ell-1} \\ 0_{\ell \times \ell-1} & U_\ell & \cdots & 0_{\ell \times \ell-1} \\ \vdots & \ddots & \ddots & \vdots \\ 0_{\ell \times \ell-1} & 0_{\ell \times \ell-1} & \cdots & U_\ell \end{pmatrix},$$

with $k$ diagonal blocks, and $\widetilde{\Gamma}_k = U_\ell^{[k]\top} \Gamma_k U_\ell^{[k]}$. Then $\widetilde{\Gamma}_k$ is the matrix of covariances generated by a system like (3.1) and (3.2), but with system matrices $(\widetilde{A}, \widetilde{C})$ and covariance matrices $(\widetilde{S}, \widetilde{R})$. Indeed, using that $\mathbf{S}_n^\perp$ is invariant under both $S$ and $A$ and that $S\Pi_{\mathbf{S}_n^\perp} = S$, one finds that for $j > i$ say,

$$\widetilde{C}\widetilde{A}^{(j-i)\vee 0}\widetilde{S}(\widetilde{A}^\top)^{(i-j)\vee 0}\widetilde{C}^\top + \delta_{ij}\widetilde{R} = U_\ell^\top \Gamma_{ij} U_\ell,$$

where $\Gamma_{ij}$ is block $(i,j)$ of $\Gamma$. Likewise, one can show that $\overline{\Gamma}_k = U_\ell^{[k]} \widetilde{\Gamma}_k U_\ell^{[k]\top}$.

CONDITION C. *The eigenvalue 0 of the covariance matrix $R$ is simple.*

This condition implies that $\widetilde{\Gamma}_k$ is nonsingular for all $k$. It also implies that the spectral density matrix $\widetilde{f}(\omega)$ of the system $(\widetilde{A}, \widetilde{C}; \widetilde{S}, \widetilde{R})$, given by

$$\widetilde{f}(\omega) = \frac{1}{2\pi}\left\{\sum_{t=-\infty}^{\infty} e^{-i\omega t}\widetilde{C}\widetilde{A}^{(-t\vee 0)}\widetilde{S}(\widetilde{A}^\top)^{(t\vee 0)}\widetilde{C}^\top + \delta_{0t}\widetilde{R}\right\},$$

has eigenvalues bounded away from zero. Since $\widetilde{A}$ is stable the eigenvalues of $\widetilde{f}(\omega)$ are bounded away from infinity as well. These observations allow us to make the stronger conclusion that the eigenvalues of $\widetilde{\Gamma}_k$ are bounded away from zero and infinity uniformly in $k$ (cf. [15], pages 265–266). In fact, the same holds for the infinite-dimensional matrix $\widetilde{\Gamma}$, which is thus also invertible. Now, put $\Gamma_k^\dagger = U_\ell^{[k]}\widetilde{\Gamma}_k^{-1}U_\ell^{[k]\top}$. Using $U_\ell^\top U_\ell = I_{\ell-1}$, one readily verifies that (i) $\Gamma_k\Gamma_k^\dagger\Gamma_k = \Gamma_k$, (ii) $\Gamma_k^\dagger\Gamma_k\Gamma_k^\dagger = \Gamma_k^\dagger$, (iii) $\Gamma_k\Gamma_k^\dagger$ is symmetric and (iv) $\Gamma_k^\dagger\Gamma_k$ is symmetric (the latter two matrices are both $\Pi_{\mathbf{S}_\ell^\perp}^{[k]}$). This implies, as the notation suggests, that $\Gamma_k^\dagger$ is the Moore–Penrose inverse of $\Gamma_k$ [25], pages 167–168. Similarly, $\Gamma^\dagger = U_\ell^{[\infty]}\widetilde{\Gamma}^{-1}U_\ell^{[\infty]\top}$.

A slightly different way to view the same properties is to note that for any vector $v \in \mathbf{S}_\ell^k$, $v^\top\Gamma_k = 0$. Thus, zero is an eigenvalue of $\Gamma_k$ with multiplicity $k$. The remaining eigenvalues agree with those of $\widetilde{\Gamma}_k$, and the assumption above guarantees that none of them is zero. Since $\mathbf{1}_\ell^\top \overline{y}_t = 0$, we also find that $v^\top \overline{Y}^- = 0$, and hence $v^\top \widehat{\Gamma}_k = 0$ for any $v \in \mathbf{S}_\ell^k$. Thus, zero is an eigenvalue with multiplicity $k$ of $\widehat{\Gamma}_k$ as well, and its remaining eigenvalues are those of $\widehat{\widetilde{\Gamma}}_k = U_\ell^{[k]\top}\widehat{\Gamma}_k U_\ell^{[k]}$. These properties imply $\|\Gamma_k\| = \|\widetilde{\Gamma}_k\|$, $\|\widehat{\Gamma}_k\| = \|\widehat{\widetilde{\Gamma}}_k\|$ and $\|\widehat{\Gamma}_k - \Gamma_k\| = \|\widehat{\widetilde{\Gamma}}_k - \widetilde{\Gamma}_k\|$. Moreover, $\widehat{\Gamma}_k^\dagger = U_\ell^{[k]}\widehat{\widetilde{\Gamma}}_k^{-1}U_\ell^{[k]\top}$ if $\widehat{\widetilde{\Gamma}}_k^{-1}$ is nonsingular.

**7. Consistency.** In this section, we give a consistency result for the parameter estimates of Section 5. The outline of the proof to a large extent follows that of [10]. First, in Section 7.1, we establish the rate of convergence of the sample covariances of the centred process. Section 7.2 contains results on consistency of the linear prediction matrix estimates, which are needed to carry out the proof in Section 7.3 of consistency of the subspace algorithm. Throughout this section, we assume that Conditions A–C hold.

7.1. *Convergence rates of sample covariances.* The following result establishes convergence rates of sample covariances, uniformly over some range of lags. We will write $z_t = o(g_t)$, if $z_t/g_t \to 0$ a.s. as $t \to \infty$.

THEOREM 1. *For any $\delta > 0$ and $k_T \leq T^\alpha$, with $\alpha = r/(2(r-2))$ and $r > 4$,*

$$(7.1) \quad \max_{-k_T < s < k_T}\left\|\frac{1}{T}\sum_{t=1}^T \overline{y}_t \overline{y}_{t+s}^\top - \mathbb{E}[\overline{y}_0 \overline{y}_s^\top]\right\| = O(Q_T) \quad \text{as } T \to \infty,$$

*with* $Q_T = T^{-1/2}(k_T \log T)^{2/r}(\log \log T)^{2(1+\delta)/r}$.

PROOF. Put $\overline{x}_t^M = x_t - \pi$. Subtracting $\pi$ and $C\pi$ from (3.1) and (3.2) yields

$$\overline{x}_{t+1}^M = A\overline{x}_t^M + \xi_{t+1} = \overline{A}\overline{x}_t^M + \xi_{t+1}, \qquad \overline{y}_{t+1} = C\overline{x}_t^M + \eta_{t+1}.$$

Using this representation, the proof proceeds essentially as that of Theorem 1 in [17]. The main difference is that the process $S_\tau(j,k,t)$ appearing therein is not always a martingale, which in turn invalidates the application of Doob's and Burkholder's inequalities. Instead, one must rely on corresponding inequalities for strongly mixing processes (Theorem 2 of [11], page 26). A complete proof can be found in a preprint version of the present paper [2]. □

7.2. *Convergence of moment matrix estimates.* Theorem 1 has immediate applications to the convergence of the moment matrix estimates.

PROPOSITION 1. *Under the assumptions of Theorem 1, and provided $k = k_T$ is such that $kO(Q_T) \to 0$, it holds that $\|\widehat{\beta}_k - \beta_k\| = kO(Q_T)$.*

PROOF. We shall first prove that $\|\widehat{\mathcal{H}}_k - \mathcal{H}_k\|$ and $\|\widehat{\Gamma}_k - \Gamma_k\|$ are both $kO(Q_T)$. Block $(i,j)$ of $\widehat{\mathcal{H}}_k - \mathcal{H}_k$ is

$$\frac{1}{T}\sum_{u=1}^T (y_{i+u-1} - m_{Y+;i})(y_{-j+u} - m_{Y-;j})^\top - \mathbb{E}[\overline{y}_0 \overline{y}_{-(i+j)}^\top]$$

$$= \frac{1}{T}\sum_{u=1}^T \{\overline{y}_{i+u-1}\overline{y}_{-j+u}^\top - \mathbb{E}[\overline{y}_0\overline{y}_{-(i+j)}^\top]\} + (C\pi - m_{Y+;i})(C\pi - m_{Y-;j})^\top,$$

where $m_{Y+;i}$ is the $i$th block of $m_{Y+}$, etc. The sum on the right-hand side agrees with that of (7.1), with $s = -(i+j)+1$, except for $2(i-1)$ terms with the lowest and highest indices. This remainder term is thus $O(k/T)$, which in turn is $O(Q_T)$. The whole sum is therefore $O(Q_T)$. We note in passing that the maximal (in absolute value) $s$ considered this way is $-2k+1$, while $|s| < k$ is required in Theorem 1; however, modifying the range of $s$ by a constant multiple (here 2) does not affect the validity of the theorem.

Regarding the second term on the right-hand side, each factor is, apart from a remainder term of order $O(k/T)$ as above, $O((\log \log T/T)^{1/2})$ by a law of iterated logarithm for strongly mixing processes [23], Theorem 5. The whole second term is therefore $O(\log \log T/T) + O(k^2/T^2)$, which again is readily checked to be $O(Q_T)$. Thus, each block of $\widehat{\mathcal{H}}_k - \mathcal{H}_k$ is $O(Q_T)$ and we obtain, for example, by bounding $\|\widehat{\mathcal{H}}_k - \mathcal{H}_k\|$ by its Frobenius norm, that $\|\widehat{\mathcal{H}}_k - \mathcal{H}_k\| = kO(Q_T)$. In a similar fashion, we find $\|\widehat{\Gamma}_k - \Gamma_k\| = kO(Q_T)$.

Corollary 8.1.6 of [13], page 396, shows that $\|\lambda_i(\widehat{\Gamma}_k) - \lambda_i(\Gamma_k)\| \leq \|\widehat{\Gamma}_k - \Gamma_k\|$, where $\lambda_i(\cdot)$ is the $i$th largest singular value of a matrix. We note that $\Gamma_k$ and $\widehat{\Gamma}_k$ are positive semi-definite, and hence their respective singular values and eigenvalues coincide. As noted in the previous section, both $\Gamma_k$ and $\widehat{\Gamma}_k$ have the eigenvalue zero with multiplicity (at least) $k$. We also saw that the multiplicity is exactly $k$ for $\Gamma_k$, and that its smallest nonzero eigenvalue is bounded away from zero in $k$. If $kO(Q_T) \to 0$ as $T \to \infty$, we can thus draw the same conclusions about $\widehat{\Gamma}_k$ from the above inequality.

For $T$ large enough that the eigenvalue zero of $\widehat{\Gamma}_k$ has multiplicity $k$, write

$$\|\widehat{\Gamma}_k^\dagger - \Gamma_k^\dagger\| = \|\widehat{\Gamma}_k^\dagger \Gamma_k \Gamma_k^\dagger - \widehat{\Gamma}_k^\dagger \widehat{\Gamma}_k \Gamma_k^\dagger\| \leq \|\widehat{\Gamma}_k^\dagger\| \|\widehat{\Gamma}_k - \Gamma_k\| \|\Gamma_k^\dagger\|,$$

where the first equality is checked as in the previous section. Here, $\|\Gamma_k^\dagger\|$ equals the reciprocal of the smallest nonzero singular value of $\Gamma_k$, whence $\|\Gamma_k^\dagger\|$ and $\|\widehat{\Gamma}_k^\dagger\|$ are both bounded as $T \to \infty$, and the middle factor is $kO(Q_T)$. We thus obtain

$$\|\widehat{\beta}_k - \beta_k\| = \|\widehat{\mathcal{H}}_k \widehat{\Gamma}_k^\dagger - \mathcal{H}_k \Gamma_k^\dagger\| \leq \|\widehat{\mathcal{H}}_k - \mathcal{H}_k\| \|\widehat{\Gamma}_k^\dagger\| + \|\mathcal{H}_k\| \|\widehat{\Gamma}_k^\dagger - \Gamma_k^\dagger\|$$
$$= kO(Q_T)O(1) + O(1)kO(Q_T) = kO(Q_T);$$

$\|\mathcal{H}_k\|$ is bounded in $k$ since the norm of its $(i,j)$ block $\mathbb{E}[\overline{y}_{t+i}\overline{y}_{t-j+1}^\top]$ tends to zero geometrically fast as $i+j \to \infty$. The proof is complete. $\square$

The next result tells us how well $\beta_k$ approximates $\mathcal{O}_k\mathcal{K}_k$.

PROPOSITION 2. *There is a $\rho \in (0,1)$ such that $\|\beta_k - \mathcal{O}_k\mathcal{K}_k\| = O(\rho^k)$.*

PROOF. Put $\mathcal{K}_k = (K, (A-KC)K, (A-KC)^2K, \ldots, (A-KC)^{k-1}K)$, put $G_k = \mathcal{O}_k\mathcal{K}_k$ and $G = \mathcal{O}\mathcal{K}$; notice that $G = \beta$. Throughout the proof, we will use matrix subindices in parentheses to denote a "block index," where the block size is $\ell$: $I_{(k)}$ is the identity matrix with $k$ blocks (that is, of size $k\ell \times k\ell$), $\Gamma_{(1:k,k+1:\infty)}$ is the submatrix of $\Gamma$ formed by block rows of indices 1 through $k$ and block columns of indices $k+1$ and upwards, etc.

Notice that, since $\beta = \mathcal{H}\Gamma^\dagger$,

$$[\beta_k - G_k, 0_{(k\times\infty)}] - [0_{(k\times k)}, G_{(1:k,k+1:\infty)}]$$
$$= [\beta_k - G_k, -G_{(1:k,k+1:\infty)}] = [\beta_k, 0_{(k\times\infty)}] - G_{(1:k,1:\infty)}$$
$$= [\beta_k, 0_{(k\times\infty)}] - [I_{(k)}, 0_{(k\times\infty)}]G = [\beta_k, 0_{(k\times\infty)}] - [I_{(k)}, 0_{(k\times\infty)}]\mathcal{H}\Gamma^\dagger$$
$$= [\mathcal{H}_k\Gamma_k^\dagger, 0_{(k\times\infty)}] - \mathcal{H}_{(1:k,1:\infty)}\Gamma^\dagger.$$

Post-multiplying the right-hand side by $\Gamma_{(1:\infty,1:k)} = \Gamma[I_{(k)}, 0_{(k\times\infty)}]^\top$ yields

$$\mathcal{H}_k\Gamma_k^\dagger\Gamma_k - \mathcal{H}_{(1:k,1:\infty)}\Gamma^\dagger\Gamma[I_{(k)}, 0_{(k\times\infty)}]^\top = \mathcal{H}_k - \mathcal{H}_{(1:k,1:\infty)}[I_{(k)}, 0_{(k\times\infty)}]^\top$$
$$= \mathcal{H}_k - \mathcal{H}_k = 0_{(k\times k)};$$

these equalities are true since, as follows from Section 6, $\Gamma_k^\dagger$ for instance works as a proper inverse on the space $(\mathbf{S}_\ell^\perp)^k$ on which $\Gamma$ and $\Gamma_k^\dagger$ operate. We thus find the following sequence of equalities:

$$[\beta_k - G_k, 0_{(k\times\infty)}]\Gamma_{(1:\infty,1:k)} = [0_{(k\times k)}, G_{(1:k,k+1:\infty)}]\Gamma_{(1:\infty,1:k)},$$

$$(\beta_k - G_k)\Gamma_k = G_{(1:k,k+1:\infty)}\Gamma_{(k+1:\infty,1:k)},$$

$$\beta_k - G_k = \mathcal{O}_k \mathcal{K}_{(k+1:\infty)} \Gamma_{(k+1:\infty,1:k)} \Gamma_k^\dagger.$$

The squared norm of $\mathcal{O}_k$ equals the largest eigenvalue of the $n \times n$ matrix $\mathcal{O}_k^\top \mathcal{O}_k = \sum_{j=0}^{k-1}(\overline{A}^i)^\top C^\top C \overline{A}^i$. Because $\overline{A}$ is stable this sum is convergent, and we conclude that $\|\mathcal{O}_k\|$ is bounded in $k$. Moreover, $\|\Gamma_k^\dagger\|$ is bounded in $k$, and $\|\Gamma_{(k+1:\infty,1:k)}\|$ is bounded in $k$ since $\overline{A}^i S \overline{A}^{j\top}$ tends to zero geometrically fast as $i \wedge j \to \infty$. Recalling (Section 6) that the blocks of $\mathcal{K}$ tend to zero geometrically fast, we conclude that so does $\|\mathcal{K}_{(k+1:\infty)}\|$. Hence, $\|\beta_k - G_k\| = O(\rho^k)$ for some $\rho \in (0,1)$, and the proof is complete. $\square$

### 7.3. Consistency of parameter estimates.

THEOREM 2. *For the algorithm of Section 5 with truncation index $k = k_T$ as in Theorem 1 and satisfying $k_T \to \infty$ and $k_T Q_T \to 0$ as $T \to \infty$, there are nonsingular random matrices $\{S_T\}$, with $\|S_T\|$ and $\|S_T^{-1}\|$ bounded a.s., such that $\|\widehat{A} - S_T A S_T^{-1}\| = o(1)$, $\|\widehat{C} - C S_T^{-1}\| = o(1)$ and $\|\widehat{K} - S_T K\| = o(1)$.*

PROOF. The proof consists of several parts. First, we show that $\widehat{\mathcal{O}}_k \widehat{\mathcal{K}}_k$ is close to $\mathcal{O}_k \mathcal{K}_k$. Then we show that the estimated centred state sequence $\{\widehat{\overline{x}}_t\}$ is, up to a change of basis, close to the true centered (predicted) state sequence $\{\overline{x}_t\}$. Then we show a similar statement about the noncentered estimated and true state sequences, and finally, we show that the estimates $\widehat{A}$, $\widehat{C}$ and $\widehat{K}$ are consistent in the sense stated above.

CLAIM 1. $\|\widehat{\mathcal{O}}_k \widehat{\mathcal{K}}_k - \mathcal{O}_k \mathcal{K}_k\| = kO(Q_T) + O(\rho^{-k})$ *for some $\rho \in (0,1)$.*

To prove this claim, note that $\widehat{\mathcal{O}}_k \widehat{\mathcal{K}}_k = \widehat{\beta}_k - \widehat{U}_2 \widehat{\Lambda}_{22} \widehat{V}_2^\top$ and thus

$$\|\widehat{\mathcal{O}}_k \widehat{\mathcal{K}}_k - \mathcal{O}_k \mathcal{K}_k\| \leq \|\widehat{\beta}_k - \mathcal{O}_k \mathcal{K}_k\| + \|\widehat{U}_2 \widehat{\Lambda}_{22} \widehat{V}_2^\top\|.$$

The norm of $\widehat{U}_2 \widehat{\Lambda}_{22} \widehat{V}_2^\top$ equals its largest singular value, which is the $n$th singular value of $\widehat{U}\widehat{\Lambda}\widehat{V}^\top = \widehat{\beta}_k$. Since $\mathcal{O}_k \mathcal{K}_k$ has rank at most $n-1$, its $n$th singular value is zero, whence, again using Corollary 8.1.6 of [13], page 396,

$$\|\widehat{U}_2 \widehat{\Lambda}_{22} \widehat{V}_2^\top\| = \sigma_n(\widehat{\beta}_k) = |\sigma_n(\widehat{\beta}_k) - \sigma_n(\mathcal{O}_k \mathcal{K}_k)| \leq \|\widehat{\beta}_k - \mathcal{O}_k \mathcal{K}_k\|.$$

We conclude that

$$\|\widehat{\mathcal{O}}_k\widehat{\mathcal{K}}_k - \mathcal{O}_k\mathcal{K}_k\| \le 2\|\widehat{\beta}_k - \mathcal{O}_k\mathcal{K}_k\| \le 2\|\widehat{\beta}_k - \beta_k\| + 2\|\mathcal{O}_k\mathcal{K}_k - \beta_k\|,$$

and Claim 1 follows by Propositions 1 and 2.

CLAIM 2. *With $\overline{S}_T = \widehat{\mathcal{O}}_k^\dagger \mathcal{O}_k \Pi_{\mathbf{S}_n^\perp}$, we have $\|T^{-1}\sum_{t=1}^T (\overline{S}_T\overline{x}_t - \widehat{\overline{x}}_t)(\overline{S}_T\overline{x}_t - \widehat{\overline{x}}_t)^\top\| = o(1)$.*

First recall that $\Pi_{\mathbf{S}_n^\perp}$ is the projection matrix onto $\mathbf{S}_n^\perp$. Since $\overline{x}_t \in \mathbf{S}_n^\perp$ we thus can and will, within the proof of this claim, take $\overline{S}_T = \widehat{\mathcal{O}}_k^\dagger \mathcal{O}_k$.

In order to prove the claim, notice that since $\widehat{\mathcal{O}}_k$ has orthonormal columns, $\widehat{\mathcal{O}}_k^\dagger = \widehat{\mathcal{O}}_k^\top$, $\widehat{\mathcal{O}}_k^\dagger \widehat{\mathcal{O}}_k = I_{n-1}$, and hence $\overline{S}_T\overline{x}_t - \widehat{\overline{x}}_t = \widehat{\mathcal{O}}_k^\dagger(\mathcal{O}_k\overline{x}_t - \widehat{\mathcal{O}}_k\widehat{\overline{x}}_t)$. Moreover,

$$\mathcal{O}_k\overline{x}_t - \widehat{\mathcal{O}}_k\widehat{\overline{x}}_t = \mathcal{O}_k[\mathcal{K}_k, (A-KC)^k\mathcal{K}]\overline{y}_t^- - \widehat{\mathcal{O}}_k\widehat{\mathcal{K}}_k(y_t^-(k) - \mathbf{1}_k \otimes m_Y)$$
$$= (\mathcal{O}_k\mathcal{K}_k - \widehat{\mathcal{O}}_k\widehat{\mathcal{K}}_k)\overline{y}_t^-(k) + \mathcal{O}_k(A-KC)^k\mathcal{K}\overline{y}_{t-k}^-$$
$$+ \widehat{\mathcal{O}}_k\widehat{\mathcal{K}}_k(\mathbf{1}_k \otimes m_Y - \mathbf{1}_k \otimes (C\pi)) = z_t^1 + z_t^2 + z_t^3$$

say. Using the orthonormality of $\widehat{\mathcal{O}}_k$, the norm of interest is bounded by

$$\sum_{i=1}^3 \left\|\sum_{t=0}^{T-1} z_t^i z_t^{i\top}\right\| + 2\sum_{i<j}\left\|\sum_{t=0}^{T-1} z_t^i z_t^{i\top}\right\|^{1/2}\left\|\sum_{t=0}^{T-1} z_t^j z_t^{j\top}\right\|^{1/2},$$

where we used the Cauchy–Schwarz inequality for the cross products. First, $\|T^{-1}\sum_{t=0}^{T-1} z_t^1 z_t^{1\top}\|$ is bounded by

$$\|\mathcal{O}_k\mathcal{K}_k - \widehat{\mathcal{O}}_k\widehat{\mathcal{K}}_k\|^2 \|\check{\Gamma}_k\| = [(kO(Q_T))^2 + O(\rho^{-2k})]O(1) = o(1)$$

by Claim 1 and the assumptions of the theorem. Here, $\check{\Gamma}_k$ is a sample covariance matrix similar to $\widehat{\Gamma}_k$ but with observations centred using the true (unknown) mean; $\|\check{\Gamma}_k\|$ is bounded in $k$ because so is $\|\Gamma_k\|$ (see Section 6) and, similar to the proof of Proposition 1, $\|\check{\Gamma}_k - \Gamma_k\| = o(1)$.

Next, we recall that $(A - KC)^j K = J^j K$ where $J$ is stable. Thus, $z_t^2 = \mathcal{O}_k J^k [K, JK, J^2 K, \ldots] y_{t-k}^-$, and hence both $\|z_t^2\|$ and $\|T^{-1}\sum_1^T z_t^2 z_t^{2\top}\|$ tend to zero geometrically fast in $k$.

Finally, $z_t^3$ does not depend on $t$, whence $\|T^{-1}\sum_{t=0}^{T-1} z_t^3 z_t^{3\top}\|$ is bounded by

$$\|\widehat{\mathcal{O}}_k\widehat{\mathcal{K}}_k\|^2 \|\mathbf{1}_k \otimes m_Y - \mathbf{1}_k \otimes (C\pi)\|^2 = \|\widehat{\mathcal{O}}_k\widehat{\mathcal{K}}_k\|^2 k\|m_Y - C\pi\|^2.$$

Here, Claim 1 shows that the first factor is bounded and, as in the proof of Proposition 1, a law of iterated logarithm ensures that the last factor is $O(\log \log T/T)$. It can be checked that $kO(\log \log T/T) = o(1)$, so it follows that the left-hand side tends to zero.

The proof of Claim 2 is complete. In fact, the argument contains a small error, in that $\widehat{\overline{x}}_t$ does not equal $\widehat{\mathcal{K}}_k(y_k^-(t) - \mathbf{1}_k \otimes m_Y)$ for $t < k$, as $y_s$ for $s \leq 0$ are not observed. Adjusting for this error gives extra sums with $k$ terms, and thus remainder terms of order $O(k/T)$. These are however of smaller-order than the main terms above.

CLAIM 3. *With*
$$S_T = \begin{bmatrix} \overline{S}_T & \mathbf{0}_{n-1} \\ \mathbf{0}_n^\top & 1 \end{bmatrix} \begin{bmatrix} I_n - \pi \mathbf{1}_n^\top \\ \mathbf{1}_n^\top \end{bmatrix},$$
*it holds that* $\|T^{-1} \sum_{t=1}^T (S_T x_t - \widehat{x}_t)(S_T x_t - \widehat{x}_t)^\top\| = o(1)$.

This claim follows easily upon observing that $S_T$ is a composition of two mappings acting as follows. Since $\mathbf{1}_n^\top x_t = 1$, the first matrix maps $x_t$ into a vector in $\mathbf{R}^{n+1}$ whose first $n$ components are $\overline{x}_t = x_t - \pi$, and whose $(n+1)$st component is 1. The second mapping maps the first $n$ components $\overline{x}_t$ into $\overline{S}_T \overline{x}_t \in \mathbf{R}^{n-1}$ and adjoins the last component 1. The output of the concatenated mappings is thus $[(\overline{S}_T \overline{x}_t)^\top, 1]^\top$. Since the last component of $\widehat{\overline{x}}_t$ is 1 by construction, the difference $S_T x_t - \widehat{x}_t$ equals $\overline{S}_T \overline{x}_t - \widehat{\overline{x}}_t$, with a final component 0 adjoined. The claim now follows from Claim 2.

CLAIM 4. *The $(n-1) \times n$ matrices $\overline{S}_T$ are such that their singular values are bounded away from infinity and zero as $T \to \infty$.*

First, we note that, since $\widehat{\mathcal{O}}_k$ has orthonormal columns, $\|\widehat{\mathcal{O}}_k\| = 1$. Likewise, $\|\Pi_{\mathbf{S}_n^\perp}\| = 1$, so that $\|\overline{S}_T\| \leq \|\mathcal{O}_k\|$. Here, $\|\mathcal{O}_k\|$ is bounded in $k$ (see the proof of Proposition 2), whence $\|\overline{S}_T\|$ is bounded from above.

We must now show that the $(n-1)$st singular value of $\overline{S}_T$ is bounded away from zero as $T \to \infty$. Following [10], page 1870, we first note that since the dimension of $\mathcal{O}_k \mathcal{K}_k$ grows with $k$, it holds that $\lambda_{n-1}(\mathcal{O}_{k+1} \mathcal{K}_{k+1}) \geq \lambda_{n-1}(\mathcal{O}_k \mathcal{K}_k)$ for all $k$. Thus, because $\mathcal{OK}$ has rank $n-1$, $\lambda_{n-1}(\mathcal{O}_k \mathcal{K}_k)$ is bounded away from zero for $k$ sufficiently large. Using that $\lambda_{n-1}(\widehat{\beta}_k) \geq \lambda_{n-1}(\mathcal{O}_k \mathcal{K}_k) - \|\widehat{\beta}_k - \mathcal{O}_k \mathcal{K}_k\|$ (Corollary 8.1.6 of [13], page 396), we find by invoking Propositions 1 and 2 that $\lambda_{n-1}(\widehat{\beta}_k)$ is bounded away from zero as well. Now $\lambda_{n-1}(\widehat{\beta}_k) = \lambda_{n-1}(\widehat{\Lambda}_{11})$, and since $\widehat{V}_1$ has orthonormal columns, the singular values of $\widehat{\mathcal{K}}_k = \widehat{\Lambda}_{11} \widehat{V}_1^\top$ and $\widehat{\Lambda}_{11}$ agree. Thus, $\lambda_{n-1}(\widehat{\mathcal{K}}_k)$ is bounded away from zero.

To complete the proof of Claim 4, let $u^\top$ be the $(n-1)$st left singular vector of $\overline{S}_T = \widehat{\mathcal{O}}_k^\dagger \mathcal{O}_k$ with unit length, that is, $u^\top \overline{S}_T = \lambda_{n-1}(\overline{S}_T) u^\top$. Then
$$\lambda_{n-1}(\overline{S}_T) \|\mathcal{K}_k \widehat{\mathcal{K}}_k^\dagger\| \geq |u \widehat{\mathcal{O}}_k^\dagger \mathcal{O}_k \Pi_{\mathbf{S}_n^\perp} \mathcal{K}_k \widehat{\mathcal{K}}_k^\dagger| = |u + u \widehat{\mathcal{O}}_k^\dagger E \widehat{\mathcal{K}}_k^\dagger|$$
$$\geq |u| - |u \widehat{\mathcal{O}}_k^\dagger E \widehat{\mathcal{K}}_k^\dagger| \geq 1 - \|\widehat{\mathcal{O}}_k^\dagger\| \|E\| \|\widehat{\mathcal{K}}_k^\dagger\|,$$

where $E = \mathcal{O}_k \Pi_{\mathbf{S}_n^\perp} \mathcal{K}_k - \widehat{\mathcal{O}}_k \widehat{\mathcal{K}}_k$. In addition $\Pi_{\mathbf{S}_n^\perp} \mathcal{K}_k = \mathcal{K}_k$ because $I_n = \Pi_{\mathbf{S}_n^\perp} + \Pi_{\mathbf{S}_n}$ and $\mathbf{1}_n^\top$ is a left eigenvector of $\mathcal{K}_k$ with eigenvalue 0 (cf. Section 6), so that $E = \mathcal{O}_k \mathcal{K}_k - \widehat{\mathcal{O}}_k \widehat{\mathcal{K}}_k$. Since $\|\widehat{\mathcal{O}}_k^\dagger\| = 1$, $\|\widehat{\mathcal{K}}_k^\dagger\|$ is bounded and $\|E\| \to 0$, we find that $\lambda_{n-1}(\overline{S}_T)$ is bounded away from zero and Claim 4 is proved.

CLAIM 5. *The $n \times n$ matrices $S_T$ are such that their singular values are bounded away from infinity and zero as $T \to \infty$.*

The squared singular values of $S_T$ equal the eigenvalues of the $n \times n$ matrix
$$S_T^\top S_T = (I_n - \mathbf{1}_n \pi^\top) \overline{S}_T^\top \overline{S}_T (I_n - \pi \mathbf{1}_n^\top) + \mathbf{1}_n \mathbf{1}_n^\top.$$
The smallest eigenvalue of this matrix is the minimum of $x^\top (S_T^\top S_T) x$ over $x \in \mathbf{R}^n$ such that $|x| = 1$. Pick $x \in \mathbf{R}^n$ with $|x| = 1$, and represent this vector (uniquely) as $x = x_1 + x_2$, where $x_1 \in \mathbf{S}_n$ and $x_2 \in \mathbf{S}_n^\perp$; also, $|x_1|^2 + |x_2|^2 = 1$. Furthermore, write $x_1 = as$ where $s = \mathbf{1}_n/\sqrt{n}$, so that $|x_1| = a$. We find that $(I_n - \pi \mathbf{1}_n^\top) x_1 = aw$, where $w = (I_n - \pi \mathbf{1}_n^\top) s \in \mathbf{S}_n^\perp$. Moreover, $(I_n - \pi \mathbf{1}_n^\top) x_2 = x_2$. Thus,

$$(7.2) \qquad x^\top (S_T^\top S_T) x = a^2 s^\top \mathbf{1}_n \mathbf{1}_n^\top s + (aw + x_2)^\top (\overline{S}_T^\top \overline{S}_T)(aw + x_2).$$

Consider $\overline{S}_T$ as a linear mapping from $\mathbf{S}_n^\perp$ to $\mathbf{R}^{n-1}$. Then there is a linear mapping $\overline{S}_T^*$ from $\mathbf{R}^{n-1}$ to $\mathbf{S}_n^\perp$, called the *adjoint* of $\overline{S}_T$, which satisfies and is defined by the relation $\langle \overline{S}_T u, v \rangle = \langle u, \overline{S}_T^* v \rangle$ for all $u \in \mathbf{S}_n^\perp$ and $v \in \mathbf{R}^{n-1}$ (e.g., [14], page 216); here $\langle \cdot, \cdot \rangle$ is the inner product of a linear space. It is straightforward to check that $\overline{S}_T^*$ is represented by the matrix $\overline{S}_T^\top$, and thus $\overline{S}_T^\top \overline{S}_T$ represents the quadratic form $\overline{S}_T^* \overline{S}_T$ acting on $\mathbf{S}_n^\perp$. This quadratic form has $n-1$ eigenvalues given by the squared $n-1$ singular values of the matrix $\overline{S}_T$ mapping $\mathbf{S}_n^\perp$ to $\mathbf{R}^{n-1}$. By Claim 4, we know that these singular values are uniformly (over $T$) bounded from below by some $\lambda > 0$, so that the eigenvalues of $\overline{S}_T^* \overline{S}_T$ are bounded from below by $\lambda^2$. We conclude that (7.2) is bounded from below by

$$(7.3) \qquad na^2 + \lambda^2 |aw + x_2|^2 \geq (n \wedge \lambda^2)(a^2 + |aw + x_2|^2),$$

where $n \wedge \lambda^2 > 0$. Here, the factor $a^2 + |aw + x_2|^2$ is bounded away from zero because either $a = |x_1|$ is not close to zero, or otherwise $x_2 \in \mathbf{S}_n^\perp$, of squared length $|x_2| = 1 - a^2$, cannot closely approximate a small vector $aw \in \mathbf{S}_n^\perp$. We now make this argument precise. If $w = 0$, the right-hand side of (7.3), apart from the factor $n \wedge \lambda^2$, equals $|x_1|^2 + |x_2|^2 = |x|^2 = 1 > 0$. Now, assume $w \neq 0$. Then for small $a = |x_1|$, say small enough that $a^2 |w|^2 \leq (1 - a^2)/16 = |x_2|^2/16$, the right-hand side of (7.3), again apart from the first factor, is bounded from below by

$$a^2 + a^2 |w|^2 + |x_2|^2 - 2a|w||x_2| \geq a^2 + 1 - a^2 - 2 \frac{1 - a^2}{4} \geq \frac{1}{2} > 0.$$

If $a$ is larger and hence satisfies the reverse inequality, i.e., $a^2 \geq 1/(16|w|^2 + 1) > 0$, the right-hand side of (7.3) is, again apart from the first factor, bounded from below by the same number. Thus, we have shown that (7.2) is uniformly bounded from below over $|x| = 1$, and Claim 5 follows.

We finally turn to proving consistency of the estimated system matrices, up to a similarity transformation. The estimates are

$$\widehat{A} = \left(\frac{1}{T}\sum_{t=0}^{T-1} \widehat{x}_{t+1}\widehat{x}_t^\top\right)\left(\frac{1}{T}\sum_{t=0}^{T-1} \widehat{x}_t\widehat{x}_t^\top\right)^\dagger,$$

$$\widehat{C} = \left(\frac{1}{T}\sum_{t=0}^{T-1} y_{t+1}\widehat{x}_t^\top\right)\left(\frac{1}{T}\sum_{t=0}^{T-1} \widehat{x}_t\widehat{x}_t^\top\right)^\dagger,$$

$$\widehat{K} = \left(\frac{1}{T}\sum_{t=0}^{T-1} \widehat{x}_{t+1}\widehat{\varepsilon}_{t+1}^\top\right)\left(\frac{1}{T}\sum_{t=0}^{T-1} \widehat{\varepsilon}_{t+1}\widehat{\varepsilon}_{t+1}^\top\right)^\dagger,$$

with $\widehat{\varepsilon}_{t+1} = y_{t+1} - \widehat{C}\widehat{x}_t$.

Consider

$$\frac{1}{T}\sum_{t=0}^{T-1} y_{t+1}\widehat{x}_t^\top - \mathbb{E}[y_{t+1}x_t^\top]S_T^\top$$

$$= \frac{1}{T}\sum_{t=0}^{T-1} y_{t+1}(\widehat{x}_t - S_T x_t)^\top + \left(\frac{1}{T}\sum_{t=0}^{T-1} y_{t+1}x_t^\top - \mathbb{E}[y_{t+1}x_t^\top]\right)S_T^\top.$$

By Claims 3 and 5 and the ergodic theorem, the norm of the right-hand side of this display is $o(1)$. Similarly, one can show that $T^{-1}\sum_{t=1}^T \widehat{x}_t\widehat{x}_t^\top - S_T\mathbb{E}[x_tx_t^\top]S_T^\top = o(1)$ and $T^{-1}\sum_{t=1}^T \widehat{x}_{t+1}\widehat{x}_t^\top - S_T\mathbb{E}[x_{t+1}x_t^\top]S_T^\top = o(1)$. Thus

$$\widehat{C} = (\mathbb{E}[y_{t+1}x_t^\top]S_T^\top + o(1))(S_T\mathbb{E}[x_tx_t^\top]S_T^\top + o(1))^\dagger$$
$$= \mathbb{E}[y_{t+1}x_t^\top]\mathbb{E}[x_tx_t^\top]^{-1}S_T^{-1} + o(1) = CS_T^{-1} + o(1),$$

where the last equality follows from (3.3) and (3.4). With completely corresponding derivations, we can also prove the consistency of $\widehat{A}$.

For analyzing $\widehat{K}$, we first need the result $\|T^{-1}\sum_{t=1}^T (\varepsilon_{t+1} - \widehat{\varepsilon}_{t+1})(\varepsilon_{t+1} - \widehat{\varepsilon}_{t+1})^\top\| = o(1)$, which can be proved using the decomposition $\varepsilon_{t+1} - \widehat{\varepsilon}_{t+1} = (Cx_t - \widehat{C}S_T x_t) - (\widehat{C}\widehat{x}_t - \widehat{C}S_T x_t)$ and then proceeding as in Claim 2, using this claim as well as the consistency of $\widehat{C}$. One can then establish the consistency of $\widehat{K}$ similarly as for $\widehat{C}$ and $\widehat{A}$, given the following observations. First, writing $\check{x}_t = x_t^M - x_t$ for the linear prediction error of $x_t^M$, which is uncorrelated with $x_t$, this together with (3.4) yields $\mathbb{E}[\varepsilon_{t+1}\varepsilon_{t+1}^\top] = R + CVC^\top$, where $V$ is the covariance matrix of $\check{x}_t$ (cf. [3]). The matrix $V$ has, just as $R$, an eigenvalue 0 with eigenvector $\mathbf{1}_\ell$. Thus, so has $\mathbb{E}[\varepsilon_{t+1}\varepsilon_{t+1}^\top]$ as well, and under

Condition C this eigenvalue is simple. Second, because $K = AVC^\top(R + CVC^\top)^\dagger$ ([3], equation (13)), the rows of $K$ are in $\mathbf{S}_\ell^\perp$. Combining these two facts, we find that

$$\mathbb{E}[x_{t+1}\varepsilon_{t+1}^\top](\mathbb{E}[\varepsilon_{t+1}\varepsilon_{t+1}^\top])^\dagger = K\mathbb{E}[\varepsilon_{t+1}\varepsilon_{t+1}^\top](\mathbb{E}[\varepsilon_{t+1}\varepsilon_{t+1}^\top])^\dagger = K,$$

because the product $\mathbb{E}[\varepsilon_{t+1}\varepsilon_{t+1}^\top](\mathbb{E}[\varepsilon_{t+1}\varepsilon_{t+1}^\top])^\dagger$ is the identity as a mapping on $\mathbf{S}_\ell^\perp$ (cf. [3], Lemma 1). Third, the row vector $\mathbf{1}_\ell^\top \widehat{C}$ estimates the coefficients of the optimal linear prediction of $\mathbf{1}_\ell^\top y_{t+1}$ from $\widehat{x}_t$. Since $\mathbf{1}_\ell^\top y_{t+1} = 1$ and the last coordinate of $\widehat{x}_t$ is 1, we have $\mathbf{1}_\ell^\top \widehat{C} = [\mathbf{0}_{n-1}^\top 1]$. Thus, $\mathbf{1}_\ell^\top \widehat{\varepsilon}_{t+1} = \mathbf{1}_\ell^\top(y_{t+1} - \widehat{C}\widehat{x}_t) = 0$, i.e. $\widehat{\varepsilon}_{t+1} \in \mathbf{S}_\ell^\top$. Therefore, the convergence

$$\frac{1}{T}\sum_{t=0}^{T-1} \widehat{\varepsilon}_{t+1}\widehat{\varepsilon}_{t+1}^\top - \frac{1}{T}\sum_{t=0}^{T-1}\varepsilon_{t+1}\varepsilon_{t+1}^\top = o(1)$$

implies that the difference of the corresponding pseudo-inverses is $o(1)$, as both sums have a simple eigenvalue 0 with the same eigenvector $\mathbf{1}_\ell$. $\square$

Our next result shows that linear predictors built from the matrices $\widehat{A}$, $\widehat{C}$ and $\widehat{K}$, are consistent without any change of basis, in the sense that they converge to optimal linear predictors given the true system matrices. Consider an observed history $\mathbf{y} = (y_0, y_{-1}, y_{-2}, \ldots)$. It follows from (4.1) and (4.2) that the optimal linear predictor of $y_m$, that is, the function $\phi_m^{\mathrm{lin}}(\mathbf{y})$ minimising $\mathbb{E}|y_m - h(\mathbf{y})|^2$ over all linear functions $h$ of $\mathbf{y}$, is

$$(7.4) \qquad \phi_m^{\mathrm{lin}}(\mathbf{y}) = CA^{m-1}\sum_{j=0}^\infty (A - KC)^j(y_{-j} - \mathbb{E}[y_0]) + \mathbb{E}[y_0].$$

In practice, one of course computes the sum above, which equals $\overline{x}_0$, recursively as additional observations become available. Since $y_m$ has binary elements, $\phi_m^{\mathrm{lin}}(\mathbf{y})$ approximates the vector of conditional probabilities $\mathbb{P}(y_m = \cdot \mid \mathbf{y})$. These conditional probabilities, which govern the optimal predictor $\phi_m^{\mathrm{opt}}(\mathbf{y})$ minimizing $\mathbb{E}|y_m - h(\mathbf{y})|^2$ over all functions $h$ of $\mathbf{y}$, are computed as $CA^{m-1}\mathbb{P}(x_1^M = \cdot \mid \mathbf{y})$ where the conditional probabilities $\mathbb{P}(x_1^M = \cdot \mid \mathbf{y})$ are essentially given by the filter (the forward pass of the forward–backward algorithm for HMMs); see, for example, [8], page 56.

THEOREM 3. *Let $\mathbf{z} = (z_0, z_{-1}, z_{-2}, \ldots)$ be a sequence in $\{d_1, d_2, \ldots, d_\ell\}$, and consider the mapping*

$$(7.5) \qquad \widehat{\phi}_m^{\mathrm{lin}}(\mathbf{z}) = \widehat{C}\widehat{A}^{m-1}\sum_{j=0}^\infty (\widehat{A} - \widehat{K}\widehat{C})^j(z_{-j} - \mathbb{E}[y_0]) + \mathbb{E}[y_0].$$

*Then, given the assumptions of Theorem 2, this mapping converges a.s. to the optimal linear $m$-step predictor for the system $(A, C)$, in the sense that (7.5) converges a.s. to (7.4), with $z_j$ instead of $y_j$, as $T \to \infty$, uniformly over all sequences $\mathbf{z}$.*

Before giving the proof, we remark that we can replace the mean $\mathbb{E}[y_0]$ by an estimate such as the sample mean $T^{-1}\sum_0^{T-1} y_t$, which converges a.s. to $\mathbb{E}[y_0]$. We also remark that since $\mathbf{1}_\ell^\top \widehat{C} = [\mathbf{0}_{n-1}^\top 1]$ (see the final part of the proof of Theorem 2) and the last row of $\widehat{A}$ (and hence of all $\widehat{A}^{m-1}$) equals $[\mathbf{0}_{n-1}^\top 1]$ (see Section 5), we have $\mathbf{1}_\ell^\top \widehat{C}\widehat{A}^{m-1} = [\mathbf{0}_{n-1}^\top 1]$. Moreover, from the proof below it follows that $\mathbf{1}_\ell^\top \widehat{C}\widehat{A}^{m-1}(\widehat{A} - \widehat{K}\widehat{C})^j \widehat{K} = 0$. We thus conclude that the prediction in (7.5) is a vector with elements summing to unity, as they should. We *conjecture*, but have not been able to prove, that the elements of this vector are in fact in $[0,1]$. Indeed, in the simulations reported in Section 8, not a single computed $\widehat{\phi}_m^{\text{lin}}(\mathbf{z})$ contained elements outside this range.

PROOF OF THEOREM 3. Since $\|z_k\|_\infty = 1$ and $(A - KC)^j K = J^j K$ where $J$ is stable (cf. the proof of Theorem 2), we find that for any $\delta > 0$ there is an $M > 0$ such that truncating the sum in (7.4) after $M$ terms yields an error less than $\delta$. Now, using Theorem 2, one finds that the sum over the first $M$ terms in (7.5) converges to the corresponding finite sum of terms in (7.4) (with $z_j$), uniformly because the $z_j$ are uniformly bounded.

What now remains to prove is that by choosing $M$ large enough, the tail sum over terms $j = M, M+1, \ldots$ in (7.5) is at most $\delta$, uniformly over $\mathbf{z}$. This follows as the product $(\widehat{A} - \widehat{K}\widehat{C})^j \widehat{K}$ has a structure similar to that of $(A - KC)^j K$. First, notice that since the last coordinate of $\widehat{x}_{t+1}$ is 1, the last row of the sum $T^{-1}\sum_0^{T-1} \widehat{x}_{t+1}\widehat{\varepsilon}_{t+1}^\top$ in the expression for $\widehat{K}$ will be the sum of the $\widehat{\varepsilon}_{t+1}^\top$. This sum is zero however, because of the definition of these variables as residuals in a linear regression. Thus, the last row of $\widehat{K}$ is zero, and so is the last row of $\widehat{K}\widehat{C}$. Moreover, the last row of $\widehat{A}$ is $[\mathbf{0}_{n-1}^\top, 1]$ (see the end of Section 5), a vector that we denote by $\widehat{\mathbf{1}}_n^\top$. Therefore, the last row of $\widehat{A} - \widehat{K}\widehat{C}$ is $\widehat{\mathbf{1}}_n^\top$. This implies that $\widehat{A} - \widehat{K}\widehat{C}$ has an eigenvalue 1 with left eigenvector $\widehat{\mathbf{1}}_n^\top$. Denote the corresponding right eigenvector by $\widehat{\gamma}$, normalized such that $\widehat{\mathbf{1}}_n^\top \widehat{\gamma} = 1$. Then $(\widehat{A} - \widehat{K}\widehat{C})^j = \widehat{J}^j + \widehat{\gamma}\widehat{\mathbf{1}}_n^\top$ say. Since $\widehat{A} - \widehat{K}\widehat{C}$ converges to $A - KC$ up to a similarity transformation, the matrix $\widehat{J}$ will be stable for large $T$, with eigenvalues converging to those of $J$. Thus, the size of the tail sum in (7.5) is of the order $\text{Sp}(\widehat{J})^M$, where Sp is the spectral radius, and we conclude that this tail sum will be less than $\delta$ for large enough $T$. As $\delta > 0$ was arbitrary, the proof is complete. □

**8. Examples.** We consider three examples, defined by the matrices

$$A_1 = \begin{pmatrix} 0.9 & 0.1 \\ 0.1 & 0.9 \end{pmatrix}, \qquad A_2 = \begin{pmatrix} 0.9 & 0.1 \\ 0.1 & 0.9 \end{pmatrix}, \qquad A_3 = \begin{pmatrix} 0.9 & 0.2 & 0.05 \\ 0.05 & 0.6 & 0.05 \\ 0.05 & 0.2 & 0.9 \end{pmatrix}$$

and

$$C_1 = \begin{pmatrix} 0.9 & 0.1 \\ 0.1 & 0.9 \end{pmatrix}, \qquad C_2 = \begin{pmatrix} 0.6 & 0.4 \\ 0.4 & 0.6 \end{pmatrix}, \qquad C_3 = \begin{pmatrix} 0.8 & 0.1 & 0.1 \\ 0.1 & 0.8 & 0.1 \\ 0.1 & 0.1 & 0.8 \end{pmatrix},$$

respectively. For $(A_1, C_1)$, the Markov chain is rather inert, and the output symbol is informative about the corresponding Markov state. For $(A_2, C_2)$, the dynamics of the hidden chain is as before, but the output symbols are much less informative about the current state. For $(A_3, C_3)$, states $e_1$ and $e_3$ are quite inert, while visits to state $e_2$ are short. The output is quite informative about the current state, although not as informative as for $C_1$.

We simulated data according to these systems, and the subspace estimation algorithm was run for increasing values of $T$ and $k$, respectively. For each system, we computed estimates $\widehat{A}$, $\widehat{C}$ and $\widehat{K}$ as in (5.1) for 250 replications of data. Then we computed the 1-step predictive distribution $\widehat{\phi}_1^{\text{lin}}$ for each time-point in an independent series of 5000 observations simulated from the same system. For each time-point $t$, we computed the $\ell_1$-norm $|\widehat{\phi}_1^{\text{lin}} - \phi_1^{\text{lin}}|_1$ comparing the estimated optimal linear predictive distribution to the true ditto, and also the $\ell_1$-norm $|\widehat{\phi}_1^{\text{lin}} - \phi_1^{\text{opt}}|_1$ comparing the same estimate to the optimal predictive distribution. These $\ell_1$-norms are also total variation distances between the respective distributions.

Table 1 reports averages of these $\ell_1$-norms over data replications and time-points for predictions. We see that the differences between $\widehat{\phi}_1^{\text{lin}}$ and $\phi_1^{\text{lin}}$ tend to zero as $T$ (and $k$) increases, which is in accordance with Theorem 3. The differences to $\phi_1^{\text{opt}}$ do not vanish in the limit however, reflecting that for each system, $\phi_1^{\text{lin}}$ and $\phi_1^{\text{opt}}$ are not equal. The limits of these differences are hence dictated by the systems $(A, C)$ themselves.

Table 2 examines the effect of the truncation index $k$, for data of length $T = 20{,}000$ from the system $(A_1, C_1)$. The subspace algorithm is reassuringly robust, in that $k$ has a negligible impact on its prediction performance, provided $k \geq 5$ say.

**9. Discussion.** The present paper is a first step towards a rigorous understanding of the use of subspace methods for estimation of and prediction with HMMs. However, many important problems remain, the most obvious one being that of obtaining a positive realization of the system matrices, that is, finding a similarity transformation that provides estimates of $A$ and $C$ with positive entries summing to unity columnwise.

The algorithm outlined in [20], Section VI.D, Step 2, based on feasible directions, finds a positive realization for a problem that is related to ours, but with a cost function that is not directly applicable to the present problem. Presumably, the algorithm could however be modified to work for the present setting. The algorithms in [9, 24] find a nonnegative matrix with a

SUBSPACE ESTIMATION METHODS FOR HMMS21

TABLE 1
*Predictive distribution errors for the $1$-step ahead linear predictor $(7.5)$, with $\mathbb{E}[y_0]$ replaced by sample means, over $250$ replications of simulated data from the systems $(A_1, C_1)$–$(A_3, C_3)$ and $5000$ out-of-sample observations. The figures are empirical averages of the $\ell_1$-norms $|\widehat{\phi}_1^{\mathrm{lin}} - \phi_1^{\mathrm{lin}}|_1$ and $|\widehat{\phi}_1^{\mathrm{lin}} - \phi_1^{\mathrm{opt}}|_1$ over the $250$ replications (for which the estimates $\widehat{A}$, etc., vary) and over the $5000$ out-of-sample observations (for which the observed history varies)*

|        |    | Mean predictive distribution error |              |              |                      |              |              |
|--------|----|------------------------------------|--------------|--------------|----------------------|--------------|--------------|
|        |    | w.r.t. optimal linear pred.        |              |              | w.r.t. optimal pred. |              |              |
| $T$    | $k$ | $(A_1, C_1)$ | $(A_2, C_2)$ | $(A_3, C_3)$ | $(A_1, C_1)$ | $(A_2, C_2)$ | $(A_3, C_3)$ |
| 1000   | 5  | 0.0299 | 0.0381 | 0.0676 | 0.0612 | 0.0384 | 0.1155 |
| 5000   | 8  | 0.0138 | 0.0173 | 0.0301 | 0.0564 | 0.0177 | 0.0992 |
| 10,000 | 12 | 0.0102 | 0.0130 | 0.0220 | 0.0562 | 0.0134 | 0.0974 |
| 20,000 | 16 | 0.0079 | 0.0095 | 0.0164 | 0.0560 | 0.0099 | 0.0964 |
| 40,000 | 20 | 0.0060 | 0.0070 | 0.0117 | 0.0559 | 0.0075 | 0.0958 |

given spectrum (set of eigenvalues), with one of them handling constraints stating that some matrix elements should be zero [9], page 1030. They do, however, not take the matrix $C$ into account, whence also these algorithms need adjustments before being suitable for the present problem. It should also be stressed that all of the above algorithms are local search algorithms, and may hence converge to a point that is not an optimal solution, or is not a solution at all. Therefore, they might need to be restarted at different initial points. That the algorithms converge only locally is not surprising as the underlying problem that they try to solve is nonconvex, but of course it removes some of the appeal of the subspace algorithm.

TABLE 2
*Results as in Table 1 for the system $(A_1, C_1)$ and various $k$*

|        |    | Mean predictive distribution error |                      |
|--------|----|------------------------------------|----------------------|
| $T$    | $k$ | w.r.t. optimal linear pred.       | w.r.t. optimal pred. |
| 20,000 | 3  | 0.0093 | 0.0550 |
| 20,000 | 5  | 0.0086 | 0.0559 |
| 20,000 | 8  | 0.0087 | 0.0559 |
| 20,000 | 10 | 0.0088 | 0.0559 |
| 20,000 | 12 | 0.0088 | 0.0559 |
| 20,000 | 15 | 0.0088 | 0.0559 |
| 20,000 | 20 | 0.0088 | 0.0559 |
| 20,000 | 25 | 0.0088 | 0.0559 |
| 20,000 | 30 | 0.0089 | 0.0559 |
| 20,000 | 40 | 0.0089 | 0.0559 |
| 20,000 | 50 | 0.0091 | 0.0558 |

Other further aspects of the subspace approach are the asymptotic distribution of the estimators, and how to minimize their asymptotic variance; the latter question also requires an understanding of subspace algorithms more general than the one presented here in terms of pre and postmultiplication of the $\widehat{\beta}_k$ by weighting matrices prior to the SVD, and factorization of the actual SVD (cf. [10]).

**Acknowledgments.** The authors thank Dr. Dietmar Bauer for helpful correspondence, and the Associate Editor and two reviewers for constructive comments that helped to improve this paper.

AstraZeneca R&D
431 83 Mölndal
Sweden

Centre for Mathematical Sciences
Lund University
Box 118
221 00 Lund
Sweden
E-mail: tobias.ryden@matstat.lu.se